\documentclass[11pt, reqno]{amsart}
\usepackage{amssymb, amsmath, amsthm, amsfonts}
\usepackage{amsrefs}
\usepackage{verbatim}
\usepackage{graphics}

\newcommand{\demo}{\begin{proof}}
\newcommand{\edemo}{\end{proof}}
\newcommand{\demoname}[1]{\begin{proof}[#1]}
\newcommand{\edemoname}{\end{proof}}
\newcommand{\demonb}{\begin{proof}}
\newcommand{\edemonb}{\renewcommand{\qedsymbol}{}\end{proof}}
\newcommand{\demonamenb}[1]{\begin{proof}[#1]}
\newcommand{\edemonamenb}{\renewcommand{\qedsymbol}{}\end{proof}}

\theoremstyle{plain}
\newtheorem{theorem}{Theorem}[section]

\newtheorem{corollary}[theorem]{Corollary}
\newtheorem{lemma}[theorem]{Lemma}

\theoremstyle{definition}
\newtheorem{example}[theorem]{Example}
\newtheorem{definition}[theorem]{Definition}
\newtheorem{remark}[theorem]{Remark}
\newtheorem{step}{Step}
\newtheorem{casee}{Case}

\newcommand{\thm}{\begin{theorem}}
\newcommand{\ethm}{\end{theorem}}

\newcommand{\expl}{\begin{example}}
\newcommand{\eexpl}{\qex\end{example}}
\newcommand{\explnb}{\begin{example}}
\newcommand{\eexplnb}{\end{example}}

\newcommand{\defn}{\begin{definition}}
\newcommand{\edefn}{\qef\end{definition}}
\newcommand{\defnnb}{\begin{definition}}
\newcommand{\edefnnb}{\end{definition}}
\newcommand{\remk}{\begin{remark}}
\newcommand{\eremk}{\qex\end{remark}}
\newcommand{\stp}{\begin{step}}
\newcommand{\estp}{\end{step}}
\newcommand{\estpp}{\qed\end{step}}
\newcommand{\cse}{\begin{casee}}
\newcommand{\ecse}{\end{casee}}

\newcommand{\coro}{\begin{corollary}}
\newcommand{\ecoro}{\end{corollary}}
\newcommand{\lem}{\begin{lemma}}
\newcommand{\elem}{\end{lemma}}

\providecommand{\qexsymbol}{$\lozenge$}%
\newcommand{\mathqex}{\quad\hbox{\qexsymbol}}
\DeclareRobustCommand{\qex}{%
  \ifmmode \mathqex
  \else
    \leavevmode\unskip\penalty9999 \hbox{}\nobreak\hfill
    \quad\hbox{\qexsymbol}%
  \fi
}
\providecommand{\qefsymbol}{$\triangle$}%
\newcommand{\mathqef}{\quad\hbox{\qefsymbol}}
\DeclareRobustCommand{\qef}{%
  \ifmmode \mathqef
  \else
    \leavevmode\unskip\penalty9999 \hbox{}\nobreak\hfill
    \quad\hbox{\qefsymbol}%
  \fi
}

\newcommand{\enum}{\begin{enumerate}}
\newcommand{\eenum}{\end{enumerate}}

\newcommand{\enumdefn}{\begin{enumerate}}
\newcommand{\eenumdefn}{\qef\end{enumerate}}

\newcommand{\nn}{\mathbb{N}}

\newcommand{\znn}{\nn \cup \{0\}}
\newcommand{\rrr}[1]{{\mathbb{R}}^{#1}}
\newcommand{\rr}{\mathbb{R}}
\newcommand{\func}[3]{{#1} \colon {#2} \to {#3}}

\newcommand{\nd}[1]{$#1$\nobreakdash-\hspace{0pt}}

\newcommand{\ncell}[1]{\nd{#1}cell}
\newcommand{\ncells}[1]{\nd{#1}cells}

\newcommand{\nface}[1]{\nd{#1}face}

\newcommand{\dsc}{discrete curvature}

\newcommand{\crc}{combinatorial Ricci curvature}
\newcommand{\Crc}{Combinatorial Ricci curvature}

\DeclareMathOperator{\ricci}{Ric}
\newcommand{\ricc}[1]{\ricci\,({#1})}

\newcommand{\ksim}[2]{{#1}^{(({#2}))}}

\newcommand{\bari}[2]{\bar {#1}_{#2}}

\newcommand{\eulgr}[1]{\chi_g ({#1})}
\newcommand{\eulp}[1]{\chi ({#1})}
\newcommand{\geulch}{ranked Euler characteristic}

\newcommand{\rkiset}[2]{{#1}_{#2}}
\newcommand{\rkinum}[2]{F_{#2}}

\newcommand{\grdd}{ranked}

\newcommand{\rfu}{rank function}

\newcommand{\setcar}[1]{|{#1}|}
\newcommand{\osetcar}[1]{|}
\newcommand{\csetcar}[1]{|}
\newcommand{\ordcom}[1]{\Delta({#1})}

\newcommand{\covered}{\prec}
\newcommand{\covers}{\succ}

\newcommand{\cfin}{covering-finite}
\newcommand{\symmdif}{\mathrel{\bigtriangleup}}

\newcommand{\rsnb}{sufficiently covered}

\newcommand{\tdap}{almost polyhedral}

\newcommand{\kd}[1]{d({#1})}
\newcommand{\kdk}[2]{\kd {{#1}, {#2}}}

\makeatletter
\@addtoreset{equation}{section}

\begin{document}

\title{Combinatorial Ricci  Curvature for Polyhedral Surfaces and Posets}
\author{Ethan D.\ Bloch}
\address{Bard College\\
Annandale-on-Hudson, NY 12504\\
U.S.A.}
\email{bloch@bard.edu}
\thanks{I would like to thank Sam Hsaio for his help with posets; and the Einstein Institute of Mathematics at the Hebrew University of Jerusalem, and especially Prof. Emanuel Farjoun, for their very kind hospitality during a sabbatical when parts of this paper were written.}
\date{}
\begin{abstract}
The \crc\ of Forman, which is defined at the edges of a CW complex, and which makes use of only the face relations of the cells in the complex, does not satisfy an analog of the Gauss-Bonnet Theorem, and does not behave analogously to smooth surfaces with respect to negative curvature.  We extend this curvature to vertices and faces in such a way that the problems with \crc\ are mostly resolved.  The discussion is stated in terms of ranked posets.
\end{abstract}
\maketitle

\markright{COMBINATORIAL RICCI CURVATURE FOR POLYHEDRAL SURFACES AND POSETS}

\section{Introduction}
\label{secINT}

There have been a number of discrete analogs of the curvature of smooth surfaces and manifolds.  The oldest such analog is the angle defect (also known as the angle deficiency) at a vertex $v$ of a triangulated polyhedral surface $M$, which is given by $\kdk vM = 2\pi - \sum_{\alpha \ni v} \alpha$, where the $\alpha$ are the angles at $v$ of the triangles containing $v$.  This curvature function goes back at least as far as Descartes (see \cite{FE}).  The angle defect satisfies various properties one would expect a curvature function to satisfy, notably an analog of the Gauss-Bonnet Theorem, which is that $\sum_{v \in M} \kdk vM = 2\pi \chi(M)$, where the summation is over all the vertices of $M$, and where $\chi(M)$ is the Euler characteristic of $M$.

The angle defect, and related constructs involving sums of angles in polyhedra, have been widely studied in dimension $2$ and higher, both for convex polytopes, for example in \cite{SH2} and \cite{GR2}, and more generally, for example in \cite{BA1}, \cite{C-M-S}, \cite{G-S2} and \cite{BL4}.

The above discrete analogs of smooth curvature are geometric in nature, making use of interior or exterior angles in simplices and polyhedral cells.  It would be interesting to know if there is a purely combinatorial definition of curvature, making use of only the face relations of the cells in the complex, that satisfies properties analogous to smooth curvature.  We consider here the combinatorial analog of curvature, called \crc, which is defined in \cite{FORM6}.  (There are other definitions of discrete curvature that make use of the term ``Ricci'' in their names, such as the ``simplicial Ricci tensor'' of \cite{A-M-M} and the ``discrete Ricci curvature'' of \cite{FRIT}, but these approaches are very different from \cite{FORM6}).  An earlier approach that is more comparable to \crc\ is the combinatorial analog of curvature given in \cite{STON1} and \cite{STON2}, and we discuss the latter briefly in Section~\ref{secDSTO}; this same approach, though with a factor of $\frac 12$, was used many years later in the context of graphs by \cite{HIGU}, and by a number of papers that refer to that one, perhaps unaware that the formula had appeared in \cite{STON1}.  

Some, though unfortunately not all, of the properties of smooth curvature have analogs for \crc.  An analog of Myers' Theorem for \crc\ is given in Theorem~6.1 of \cite{FORM6}, which says that if $K$ is an appropriately nicely behaved CW complex, and if the \crc\ is positive at every edge in $K$, then $\pi_1 (K)$ is finite; see \cite{MYER} for the original smooth version of this theorem.  Theorem~7.2 of \cite{FORM6} shows that for any simplicial complex $K$ of dimension at least $2$ that is a combinatorial manifold, there is a subdivision $M$ of $K$ such that $\ricc e < 0$ for every edge $e$ of $M$; see \cite{GAO1}, \cite{GA-YA}, \cite{BROO1}, \cite{LOHK1} and \cite{LOHK2} for the smooth analog of that result for Ricci curvature of Riemannian manifolds in dimensions $3$ or higher.  

The analogy between \crc\ and smooth curvature breaks down in dimension $2$, because there can be no smooth analog of Theorem~7.2 of \cite{FORM6} in dimension $2$, due to the Gauss-Bonnet Theorem for smooth surfaces, which implies that a smooth surface with non-negative Euler characteristic cannot have curvature that is everywhere negative.  \Crc, therefore, does not satisfy an analog of the Gauss-Bonnet Theorem, and that again shows that \crc\ is not entirely analogous to smooth curvature.  

The purpose of this note is to offer a way to resolve this anomaly of \crc\ in  dimension $2$.  Specifically, we describe a way to extend \crc\ to all cells of a \nd{2}dimensional polyhedral complex, and we  prove an analog of the Gauss-Bonnet Theorem for this extended curvature.  We also provide an answer to the question of everywhere negative curvature in the orientable case.

\Crc\ is the \nd{1}dimensional case of a more general definition of curvature in \cite{FORM6} that applies to cells in all dimensions, and which is defined as follows.  Let $\alpha$ and $\eta$ be \ncells p\ in a CW complex.  The \ncells p\ $\alpha$ and $\eta$ are called parallel neighbors if they are either both the faces of a common \ncell {(p+1)}, or they both have a common \nface {(p-1)}, but not both.  The curvature at $\alpha$ is then defined, in the notation of \cite{FORM6}, by
\begin{equation}\label{eqBBC}
\#\{\text{\ncells {(p+1)} $\beta > \alpha$}\} + \#\{\text{\ncells {(p-1)} $\gamma < \alpha$}\} - \#\{\text{parallel neighbors of $\alpha$}\},
\end{equation}
where $<$ denotes the relation of being a face.  (In \cite{FORM6} the above definition is also given with weights on the cells, though we do not do that here.)  \Crc\ is the special case of this curvature when $p = 1$, and is denoted $\ricc e$ for every edge $e$ of the CW complex. 

We will restrict our attention to \nd{2}dimensional cell complexes.  From the point of view of \crc\ and the fundamental group, restricting to dimension $2$ is no loss, because both are computed entirely in the \nd{2}skeleton of a cell complex.

Because Equation~(\ref{eqBBC}) uses only the face relations of the cells of a CW complex, it is  more clear, and slightly more general, to formulate our discussion in the context of posets, which we do in Section~\ref{secDRC}; in Section~\ref{secSUR} we will return to  \nd{2}dimensional polyhedral complexes.

Although the definition of curvature in \cite{FORM6} carries over directly to ranked posets, as long as every element covers, and is covered by, finitely many elements, the approach we take here uses a slight variant of that definition, where we replace the number of parallel neighbors with the more convenient set of all neighbors (which means that non-parallel neighbors are double counted).

\section{Discrete Curvature on Ranked Posets}
\label{secDRC}

We assume that the reader is familiar with basic properties of posets.  See \cite{STAN}*{Chapter~3} for details.  Let $P$ be a poset.  We let $<$ denote the partial order relation on $P$, and we write $a \covered b$ if $b$ covers $a$.  A function $\func {\rho}{P}{\{0, 1, \ldots, r\}}$ for some $r \in \znn$ is a \rfu\ for $P$ if it satisfies the following conditions: for $a, b \in P$, if $a$ is a minimal element then $\rho (a) = 0$, and if $a \covered b$ then $\rho (a) + 1 = \rho (b)$.    A poset is \grdd\ if it has a \rfu.  If a poset has a \rfu, the \rfu\ is unique.
The rank of such a poset is the smallest possible $r$.  If $P$ is \grdd\ and has rank $r$, and if $i \in \{0, 1, \ldots, r\}$, let $\rkiset Pi = \{x \in P \mid \rho (x) = i\}$ and $\rkinum Pi = \setcar {\rkiset Pi}$.

We note that a \grdd\ poset need not be graded, using the terminology of \cite{STAN}*{Chapter~3}; the term ``\grdd'' is used in \cite{ROS00}*{Section~11.1.3}, though we use the definition of a \rfu\ given in \cite{STAN}*{Chapter~3}.

Let $P$ be a poset.  The order complex of $P$, denoted $\ordcom P$, is the simplicial complex with a vertex for each element of $P$, and a simplex for each non-empty chain of elements of $P$.  It is a standard fact that this construction yields a simplicial complex.
Suppose that $P$ is finite.  The Euler characteristic of $P$, denoted $\eulp P$ is defined by $\eulp P = \chi(\ordcom P)$.

The reason to define $\eulp P$ as $\chi(\ordcom P)$ is that a finite poset in general has no natural \rfu, in contrast to the simplicial complex $\ordcom P$, which is naturally ranked by the dimensions of the simplices, and this natural rank function is what allows the Euler characteristic of simplicial complexes to be defined.

Suppose, however, that $P$ is a finite \grdd\ poset.  Then there is a more direct approach to defining the Euler characteristic of $P$, as given in the following definition.

\defn\label{defnADQ}
Let $P$ be a finite ranked poset of rank $r$.  The \textbf{\geulch} of $P$ is the number $\eulgr P = \sum_{i=0}^r (-1)^i \rkinum Pi$.
\edefn

If $P$ is the face poset of a finite regular CW complex, or in particular a polyhedral complex, or simplicial complex, then $\eulgr P = \eulp P$.  In general, however, it is not the case that $\eulgr P$ and $\eulp P$ are equal.

We also need the following definition regarding posets.

\defn\label{defnBCA}
Let $P$ be a poset.  The poset $P$ is \textbf{\cfin} if for any $a \in P$, there are finitely many $b \in P$ such that $a \covered b$, and there are finitely many $c \in P$ such that $c \covered a$.  
\edefn

Throughout this section, let $P$ be a \cfin\ ranked poset of rank $r$.  Let $i \in \{0, 1, \ldots, r\}$, and $x \in \rkiset Pi$.  
Let 
\begin{gather*}
A_i(x) = \setcar {\{y \in \rkiset P{i+1} \mid x \covered y\}}, \quad B_i(x) = \setcar {\{z \in \rkiset P{i-1} \mid z \covered x\}}\\
U_i(x) = \sum_{y \covers x} B_{i+1}(y), \quad D_i(x) = \sum_{z \covered x} A_{i-1}(z),
\end{gather*}
and
\begin{align*}
N_i(x) = \osetcar \{\{w \in \rkiset Pi &\mid x \covered v \text{ and } w \covered v \text{ for some } v \in \rkiset P{i+1}\}\\
 &\symmdif \{w \in \rkiset Pi \mid u \covered x \text{ and } u \covered w \text{ for some } u \in \rkiset P{i-1}\}\csetcar,
\end{align*}
where $\bigtriangleup$ denotes symmetric difference, and summation over the empty set is taken to be zero.  

The reader can verify the following equalities:
\begin{align}
\sum_{x \in \rkiset Pi} A_i(x) &= \sum_{y \in \rkiset P{i+1}} B_{i+1}(y)\label{eqBAB1}\\
\sum_{x \in \rkiset Pi} U_i(x) 
 &= \sum_{y \in \rkiset P{i+1}} [B_{i+1}(y)]^2\label{eqBAB2}\\
\sum_{x \in \rkiset Pi} D_i(x) 
 &= \sum_{z \in \rkiset P{i-1}} [A_{i-1}(z)]^2.\label{eqBAB3}
\end{align}

Equation~(\ref{eqBBC}) in the case $p = 1$, which defines \crc, can be rewritten as
\begin{equation}\label{eqRIC}
\ricc e = A_1(e) + B_1(e) - N_1(e)
\end{equation}
for all $e \in \rkiset P1$.

We now define our \dsc\ functions on \cfin\ ranked posets of rank $2$.

\defn\label{defnBAY}
let $P$ be a \cfin\ ranked poset of rank $2$.  For each $i \in \{0, 1, 2\}$, let $\func {R_i}{\rkiset Pi}{\rr}$ be defined by  
\begin{align*}
R_0(v) &= 1 + \frac 32A_0(v) - [A_0(v)]^2,
\\
R_1(e) &= 1 + 6A_1(e) + \frac 32B_1(e) - U_1(e) - D_1(e),
\\
R_2(\sigma) &= 1 + 6B_2(\sigma) - [B_2(\sigma)]^2,
\end{align*}
for all $v \in \rkiset P0$ and $e \in \rkiset P1$ and $\sigma \in \rkiset P2$.
\edefn

We will see in Lemma~\ref{lemBEBX} that for certain posets, including the face posets of all \nd{2}dimensional polyhedral complexes, the function $R_1$ equals $\ricci$.

The choice of coefficients in Definition~\ref{defnBAY}, particularly $6$ and $\frac 32$, may seem unmotivated.  They were chosen simply because they relate properly to \crc.  A slight variation in the coefficients can also be used, but the above choice appears to be the simplest possible.

The following analog of the Gauss-Bonnet Theorem is a trivial consequence of Definition~\ref{defnBAY} and Equations~(\ref{eqBAB1})--(\ref{eqBAB3}); the details are left to the reader. 

\thm\label{thmBAK}
Let $P$ be a finite ranked poset of rank $2$.  Then
\[
\sum_{v \in \rkiset P0} R_0(v) - \sum_{e \in \rkiset P1} R_1(e) + \sum_{\sigma \in \rkiset P2} R_2 (\sigma) = \eulgr P.
\]
\ethm

We now turn to a less trivial result, which says something about the nature of the poset $P$ in the case that the $R_1$ is everywhere positive, somewhat analogously to Theorem~6.1 of \cite{FORM6}.  Our result is weaker than that theorem, due to the fact that not all posets are as nicely behaved as the face posets of CW complexes.  

In the case of Gaussian curvature of compact smooth surfaces, the classical Gauss-Bonnet Theorem implies that if the curvature is everywhere positive, then the Euler characteristic of the surface is positive.  A similar result holds for the polyhedral curvature defined in \cite{BA1}, and for the combinatorial approach of \cite{STON1} and \cite{STON2} (see Section~\ref{secDSTO} for a brief discussion of that approach).  Unfortunately, because of the negative coefficient for the $R_1$ terms in Theorem~\ref{thmBAK}, it is not possible to deduce from this version of the Gauss-Bonnet Theorem that if each of $R_0$, $R_1$ and $R_2$ are everywhere positive, then the Euler characteristic is positive.  It turns out, as we now see, that it is nonetheless true that positive $R_1$ implies that the \geulch\ is positive.  In fact, all that is required is that the average value of $R_1$ is positive (in contrast to Theorem~6.1 of \cite{FORM6}, which requires $\ricci$ to be positive everywhere).  
  
We start with a definition.

\defnnb\label{defnBAV}
Let $P$ be a finite ranked poset of rank $2$.  Let
\[
\bari R1 = \frac 1{\rkinum P1} \sum_{e \in \rkiset P1} R_1(e), \qquad \bari A1 = \frac 1{\rkinum P1} \sum_{e \in \rkiset P1} A_1(e), \qquad \bari B1 = \frac 1{\rkinum P1} \sum_{e \in \rkiset P1} B_i(e).
\]
The poset $P$ is \textbf{\rsnb} if 
\begin{equation}\label{eqBEE}
[\bari A1 + \bari B1]^2 - 6\bari A1 - \frac 32\bari B1 - 1 \ge 0.\tag*{\qefsymbol}
\end{equation}
\edefnnb

\remk\label{remkCCC}
It is straightforward to verify that if $P$ is a finite ranked poset of rank $2$ and if $\bari B1 = 2$, then $P$ is \rsnb\ if and only if $\bari A1 \ge 2$;
in particular, that would hold for the face poset of any simplicial surface, and more generally for any polyhedral map of a surface (as defined in Section 21.1 of \cite{BR-SC}).  If $\bari B1 \ge \frac {20}9$, then $P$ is \rsnb\ regardless of the value of $\bari A1$.
\eremk

In the following theorem, our analog of Theorem~6.1 of \cite{FORM6}, we restrict our attention to \rsnb\ posets.  After the theorem, we will give a simple example that shows the necessity of some such restriction.  A relation between similar average values and the topology of \nd{3}manifolds is discussed in \cite{L-S}, so it is not surprising that such averages are used here as well; it is not clear whether the results of \cite{L-S} are related to our approach.  In \cite{L-S} the averages take place in simplicial complexes, so $\bari B1 = 2$, and therefore only $\bari A1$ is considered.  Moreover, while there is a very simple formula for $\bari A1$ in a simplicial complex, as used in \cite{L-S}, that is not the case for posets that are not the face posets of simplicial complexes.
 
\thm\label{thmBAL}
Let $P$ be a finite \rsnb\ ranked poset of rank $2$.  If $\bari R1 > 0$, then $\eulgr P > 0$.
\ethm

\demo
Suppose that $\bari R1 > 0$.  By Equations~(\ref{eqBAB2}) and (\ref{eqBAB3}), we see that
\begin{align*}
\frac 12 \bari R1\rkinum P1 < \bari R1\rkinum P1 &= \sum_{e \in \rkiset P1} R_1(e) = \sum_{e \in \rkiset P1} \bigl[1 + 6A_1(e) + \frac 32B_1(e) - U_1(e) - D_1(e)\bigr]\\
 &=  \rkinum P1 + 6\rkinum P1 \bari A1 + \frac 32\rkinum P1\bari B1 - \sum_{\sigma \in \rkiset P2} [B_2(\sigma)]^2 - \sum_{v \in \rkiset P0} [A_0(v)]^2,
\end{align*}
and hence
\[
\sum_{v \in \rkiset P0} [A_0(v)]^2 + \sum_{\sigma \in \rkiset P2} [B_2(\sigma)]^2 < \rkinum P1\bigl[1 - \frac 12 \bari R1 + 6\bari A1 + \frac 32\bari B1\bigr].
\]
Let $d = 1 - \frac 12 \bari R1 + 6\bari A1 + \frac 32\bari B1$.  The inequality
\[
\frac 1n \left[\sum_{i=1}^n a_i\right]^2 \le \sum_{i=1}^n \left(a_i\right)^2
\]
implies that
\[
\frac 1{\rkinum P0}\left[\sum_{v \in \rkiset P0} A_0(v)\right]^2 + \frac 1{\rkinum P2}\left[\sum_{\sigma \in \rkiset P2} B_2(\sigma)\right]^2 < \rkinum P1d.
\]
By Equation~(\ref{eqBAB1}) we deduce that 
\[
\frac 1{\rkinum P0}\left[\sum_{e \in \rkiset P1} B_1(e)\right]^2 + \frac 1{\rkinum P2}\left[\sum_{e \in \rkiset P1} A_1(e)\right]^2 < \rkinum P1d.
\]
Hence
\[
\frac 1{\rkinum P0}[\rkinum P1 \bari B1]^2 + \frac 1{\rkinum P2}[\rkinum P1 \bari A1]^2 < \rkinum P1d,
\]
and therefore
\[
[\bari B1]^2\frac {\rkinum P1}{\rkinum P0} + [\bari A1]^2\frac {\rkinum P1}{\rkinum P2} < d.
\]
Let $a = [\bari A1]^2$ and $b = [\bari B1]^2$.  Then
\begin{equation}\label{eqBAC}
b\frac {\rkinum P1}{\rkinum P0} + a\frac {\rkinum P1}{\rkinum P2} < d.
\end{equation}

Because $P$ is a ranked poset of rank $2$, it follows that $\bari A1 \ne 0$ and $\bari B1 \ne 0$.  Hence $a, b > 0$.  Because $b\frac {\rkinum P1}{\rkinum P0} + a\frac {\rkinum P1}{\rkinum P2} > 0$, then $d > 0$.

By hypothesis we know that $P$ is \rsnb\ and $\frac 12 \bari Ri > 0$, and hence
\[
\left[\bari B1 + \bari A1\right]^2 \ge 1 + 6\bari A1 + \frac 32\bari B1 > 1 - \frac 12 \bari R1 + 6\bari A1 + \frac 32\bari B1,
\]
which is the same as $(\sqrt a + \sqrt b)^2 > d$.

Let
\[
U = \{(x, y) \in \rrr{2} \mid x > 0 \text{ and } y > 0 \text{ and } ax + by < d\},
\]
and let $\func f{\rrr{2}}{\rr}$ be defined by $f(x, y) = (x - 1)(y - 1)$ for all $(x, y) \in \rrr{2}$.  The level curve for $f$ with value $c = 1$ is the hyperbola $y = \frac 1{x - 1} + 1$,
and the function $f$ has value less than $1$ between the branches of this hyperbola.
Using the condition $(\sqrt a + \sqrt b)^2 > d$,
the reader can verify that the set $U$ is between the two branches of the hyperbola (find the point on the upper branch of the hyperbola at which the tangent line is parallel to the line $ax + by = d$).  It follows that $f(x, y) < 1$ for all $(x, y) \in U$.

Combining the fact that $\frac {\rkinum P1}{\rkinum P0} > 0$ and $\frac {\rkinum P1}{\rkinum P2} > 0$ with Equation~(\ref{eqBAC}), we know that $(\frac {\rkinum P1}{\rkinum P0}, \frac {\rkinum P1}{\rkinum P2}) \in U$, and therefore $(\frac {\rkinum P1}{\rkinum P0} - 1)(\frac {\rkinum P1}{\rkinum P2} - 1) < 1$.  On the other hand, we see that
\begin{align*}
\left(\frac {\rkinum P1}{\rkinum P0} - 1\right) \left(\frac {\rkinum P1}{\rkinum P2} - 1\right) &= \frac {\rkinum P1 - \rkinum P0}{\rkinum P0} \cdot \frac {\rkinum P1 - \rkinum P2}{\rkinum P2}\\
 &= \frac {\rkinum P2 - \eulgr P}{\rkinum P0} \cdot \frac {\rkinum P0 - \eulgr P}{\rkinum P2}\\
 &= \frac 1{\rkinum P0\rkinum P2}[\eulgr P]^2 - \frac {\rkinum P0 + \rkinum P2}{\rkinum P0\rkinum P2}\eulgr P + 1.
\end{align*}
It now follows that
\[
\frac 1{\rkinum P0\rkinum P2}[\eulgr P]^2 - \frac {\rkinum P0 + \rkinum P2}{\rkinum P0\rkinum P2}\eulgr P + 1 < 1,
\]
and hence
\[
[\eulgr P]^2 - (\rkinum P0 + \rkinum P2)\eulgr P < 0.
\]
Because $x^2 - (\rkinum P0 + \rkinum P2)x < 0$ if and only if $0 < x < \rkinum P0 + \rkinum P2$, we conclude that $\eulgr P > 0$.
\edemo

\begin{figure}[ht]
\centering\includegraphics{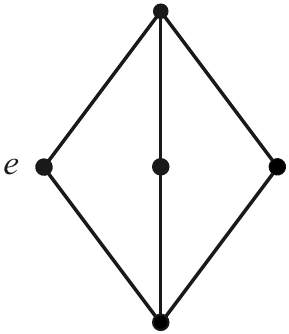}
\caption{}\label{figDIRIr}
\end{figure}

The hypothesis in Theorem~\ref{thmBAL} that $P$ is \rsnb\ cannot be dropped.  Let $P$ be the poset shown in Figure~\ref{figDIRIr}.  Then $R_1(e) = \frac 52$, and so $\bari R1 = \frac 52$, and yet $\eulgr P = -1$.  Hence, some hypothesis on $P$ is required for the theorem to hold. 

\begin{figure}[ht]
\centering\includegraphics{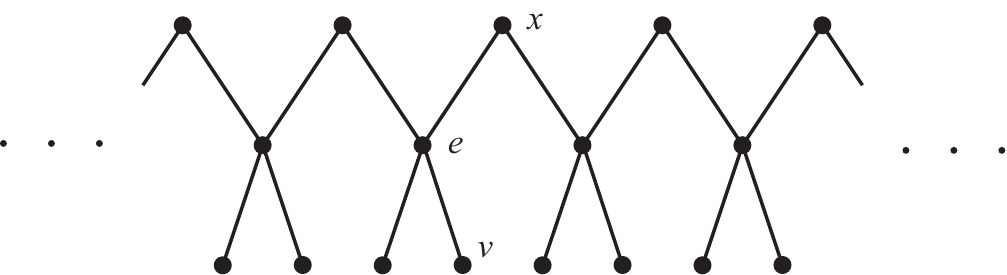}
\caption{}\label{figDIRIm}
\end{figure}

The analog of Myers' Theorem in Theorem~6.1 of \cite{FORM6} has a much stronger conclusion than our Theorem~\ref{thmBAL}.  Specifically, the main part of the proof of Theorem~6.1 of \cite{FORM6} consists of proving that if an appropriately nicely behaved connected CW complex has everywhere positive \crc\ that is bounded away from zero, then the CW complex is bounded, in the sense that there is an upper bound on the lengths of paths between vertices.  If the CW complex is regular and \cfin, that would imply that the CW complex is finite.  Unfortunately, the analogous result does not hold when $R_1$ is everywhere positive and bounded away from zero on a \cfin\ ranked poset of rank $2$, or even when all three of $R_0$, $R_1$ and $R_2$ are everywhere positive and bounded away from zero, as seen in the following example.  Let $Q$ be the poset shown in Figure~\ref{figDIRIm}, where the pattern repeats infinitely.  Let $v, e, x \in P$ be the elements shown in the figure.  Then $R_0(v) = \frac 32$, and $R_1(e) = 4$ and $R_2(x) = 9$, and yet $Q$ is infinite.  We note that $\ricci$ works no better for this poset, because $\ricc e = 2$.

Finally, we note that whereas in the special case of compact smooth surfaces, the Gauss-Bonnet Theorem immediately implies that if the average curvature is positive then so is the Euler characteristic, it is not so simple in our present context of finite ranked posets of rank $2$, because our analog of the Gauss-Bonnet Theorem, Theorem~\ref{thmBAK}, makes use of $R_0$, $R_1$ and $R_2$, whereas Theorem~\ref{thmBAL} uses only $R_1$, and hence the latter theorem does not follow from the former.

\section{$2$-Dimensional Polyhedral Complexes}
\label{secSUR}

We now relate our curvature to the problem with \crc\ that was mentioned at the start of Section~\ref{secINT}.  

First, we note that for an arbitrary \cfin\ ranked poset of rank $2$, it is not necessarily the case that $R_1$ equals $\ricci$.  For example, let $P$ be the poset shown in Figure~\ref{figDIRIr}.  Then $R_1(e) = \frac 52$, but $\ricc e = 2$.

However, as we now show, it is the case that $R_1$ equals $\ricci$ for the following class of posets.

\defnnb\label{defnBAW}
Let $P$ be a ranked poset of rank $2$.  The poset $P$ is \textbf{\tdap} if it is  \cfin, and if the following conditions hold.  Let $w \in \rkiset P0$, and $a, b \in \rkiset P1$ and $\tau \in \rkiset P2$.  Suppose $a \ne b$. 
\enumdefn
\item
$B_1(a) = 2$.
\item
There is at most one $v \in \rkiset P0$ such that $v \covered a, b$.
\item
There is at most one $\sigma \in \rkiset P2$ such that $a, b \covered \sigma$.
\item
If $w < \tau$, then $[w, \tau]$ has cardinality four.
\eenumdefn
\edefnnb

\begin{figure}[ht]
\centering\includegraphics{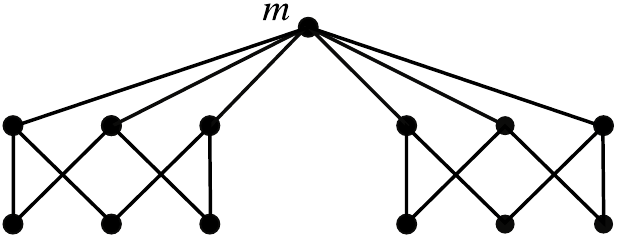}
\caption{}\label{figMCOM7}
\end{figure}

\begin{figure}[ht]
\centering\includegraphics{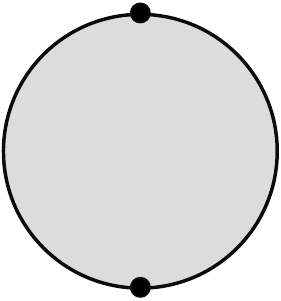}
\caption{}\label{figMCOM8}
\end{figure}

The face poset of every \nd{2}dimensional polyhedral complex, and in particular every simplicial complex, is an \tdap\ poset.  The face poset of a polyhedral map on a compact surface is also an \tdap\ poset.  However, the set of \tdap\ posets neither contains, nor is contained in, the set of face posets of all \nd{2}dimensional regular CW complexes.  The poset in Figure~\ref{figMCOM7} is \tdap\ but not the face poset of a regular CW complex, because if it were, then the boundary of $m$ would be disconnected.  On the other hand, The face poset of the CW complex with two vertices, two edges and one disk, seen in Figure~\ref{figMCOM8}, is not \tdap.  We note that Condition~(4) in Definition~\ref{defnBAW} is found in a number of places, such as Proposition~2.2 of \cite{BJOR2} and Definition~3.3 of \cite{FORM6}.

\lem\label{lemBEBX}
Let $P$ be an \tdap\ ranked poset of rank $2$.  Then $R_1(e) = \ricc e$ for all $e \in \rkiset P1$.
\elem

\demonb
Let $e \in \rkiset P1$.  Let 
\begin{align*}
U &= \{a \in \rkiset P1 \mid e \covered \tau \text{ and } a \covered \tau \text{ for some } \tau \in \rkiset P2\}\\
\intertext{and}
V &= \{a \in \rkiset P1 \mid v \covered e \text{ and } v \covered a \text{ for some } v \in \rkiset P0\}.
\end{align*}
Then $N_1(e) = \setcar {U \symmdif V}$.  

If $\sigma \in \rkiset P2$ and $e \covered \sigma$, let 
\[
C_\sigma = \{a \in \rkiset P1 \mid a \covered \sigma \text{ and } a \notin V\}
\]
and if $v \in \rkiset P0$ and $v \covered e$, let 
\[
D_v = \{a \in \rkiset P1 \mid v \covered a \text{ and } a \notin U\}.
\]
Then $U - V = \bigcup_{\sigma \covers e} C_\sigma$ and $V - U = \bigcup_{v \covered e} D_v$.  Observe that $e \in U \cap V$.  It follows that $e \notin C_\sigma$ for all $\sigma \in \rkiset P2$ such that $e \covered \sigma$, and $e \notin D_v$ for all $v \in \rkiset P0$ such that $v \covered e$.

Let $u, w \in \rkiset P0$.  Suppose that $u, w \covered e$ and $u \ne w$.  We claim that $D_u \cap D_w= \emptyset$ and that $\setcar {D_u} = A_0(u) - A_1(e) - 1$.  

Let $d \in D_u \cap D_w$.  Then $u, w \covered d$, and because $d \ne e$, we have a contradiction to Condition~(2) of Definition~\ref{defnBAW}.  Hence $D_u \cap D_w= \emptyset$.

Let $S = \{a \in \rkiset P1 \mid u \covered a\}$.  Then 
\[
S = \{a \in \rkiset P1 \mid u \covered a \text{ and } a \notin U\} \cup \{a \in \rkiset P1 \mid u \covered a \text{ and } a \in U - \{e\}\} \cup \{e\},
\]
where this union is disjoint.  

Let 
\[
\func f{\{a \in \rkiset P1 \mid u \covered a \text{ and } a \in U - \{e\}\}}{\{\sigma \in \rkiset P2 \mid e \covered \sigma\}}
\]
be defined as follows.  Let $m \in \{a \in \rkiset P1 \mid u \covered a \text{ and } a \in U - \{e\}\}$.  By the definition of $U$ there is some $\tau \in \rkiset P2$ such that $m, e \covered \tau$.  Then $\tau$ is unique by Condition~(3) of Definition~\ref{defnBAW}, and let $f(m) = \tau$.  Let $n \in \{a \in \rkiset P1 \mid u \covered a \text{ and } a \in U - \{e\}\}$, and suppose that $f(m) = f(n)$.  Then $u \covered m, e, n$ and $m, e, n \covered f(m) = f(n)$, which contradicts Condition~(4) of Definition~\ref{defnBAW}, unless $m = n$.  Hence $f$ is injective.  Let $\eta \in \{\sigma \in \rkiset P2 \mid e \covered \sigma\}$.  Then $u \covered e \covered \eta$, and by Condition~(4) of Definition~\ref{defnBAW} there is a unique $c \in \rkiset P1$ such that $c \ne e$ and $u \covered c \covered \eta$.  Then $c \in {\{a \in \rkiset P1 \mid u \covered a \text{ and } a \in U - \{e\}\}}$ and $f(c) = \eta$.  Hence $f$ is surjective.

Because $f$ is bijective, we see that
\[
\setcar S = \setcar {\{a \in \rkiset P1 \mid u \covered a \text{ and } a \notin U\}} + \setcar {\{\sigma \in \rkiset P2 \mid e \covered \sigma\}} + 1,
\]
which implies that $A_0(u) = \setcar {D_u} + A_1(e) + 1$, and hence $\setcar {D_u} = A_0(u) - A_1(e) - 1$.

Using a similar argument as above, together with the fact that $B_1(e) = 2$, which holds by Condition~(1) of Definition~\ref{defnBAW}, it is seen that if $\eta, \tau \in \rkiset P2$, and if $e \covered \eta, \tau$ and $\eta \ne \tau$, then $C_\eta \cap C_\tau = \emptyset$ and $\setcar {C_\eta} = B_2(\eta) - 3$.  

It follows that
\begin{align*}
N_1(e) &= \setcar {U - V} + \setcar {V - U} = \sum_{\sigma \covers e} \setcar {C_\sigma} + \sum_{v \covered e} \setcar {D_v}\\
 &= \sum_{\sigma \covers e} \left[B_2(\sigma) - 3\right] +\sum_{v \covered e} \left[A_0(v) - A_1(e) - 1\right].
\end{align*}

Equation~(\ref{eqRIC}) then yields
\begin{align*}
\ricc e &= A_1(e) + B_1(e) - N_1(e)\\
 &= A_1(e) + 2 - \sum_{\sigma \covers e} \left[B_2(\sigma) - 3\right] - \sum_{v \covered e} \left[A_0(v) - A_1(e) - 1\right]\\
 &= A_1(e) + 2 - \sum_{\sigma \covers e} B_2(\sigma) + 3A_1(e) - \sum_{v \covered e} A_0(v) + 2A_1(e) + 2\\
 &= 4 + 6A_1(e) - U_1(e) - D_1(e)\\
 & = 1 + \frac 32B_1(e) + 6A_1(e) - U_1(e) - D_1(e) = R_1(e).\tag*{\qedsymbol}
\end{align*}
\edemonb
\vspace{-\baselineskip}

Combining Theorem~\ref{thmBAK}, Lemma~\ref{lemBEBX} and Remark~\ref{remkCCC}, we now see that by having $R_0$ and $R_2$ available to us, there is a Gauss-Bonnet Theorem for \tdap\ ranked posets of rank $2$ that incorporates $\ricci$.

\coro\label{coroBBA}
Let $P$ be a finite \tdap\ ranked poset of rank $2$.  Then
\[
\sum_{v \in \rkiset P0} R_0(v) - \sum_{e \in \rkiset P1} \ricc e + \sum_{\sigma \in \rkiset P2} R_2 (\sigma) = \eulgr P.
\]
\ecoro

Restricting our attention to polyhedral complexes allows us to make use of the standard Euler characteristic.  If $K$ is a polyhedral complex, we let $\ksim Ki$ denote the collection of all \ncells i\ of $K$, for each $i \in \{0, 1, 2\}$.

\coro\label{coroBBY}
Let $K$ be a \nd{2}dimensional polyhedral complex.  Then
\[
\sum_{v \in \ksim K0} R_0(v) - \sum_{e \in \ksim K1} \ricc e + \sum_{\sigma \in \ksim K2} R_2 (\sigma) = \chi(P).
\]
\ecoro

The following corollary is deduced immediately from from Theorem~\ref{thmBAL}, Lemma~\ref{lemBEBX} and Remark~\ref{remkCCC}, together with the fact that $\bari B1 = 2$ in an \tdap\ ranked poset of rank $2$.

\coro\label{coroBAX}
Let $P$ be a finite \tdap\ ranked poset of rank $2$.  Suppose that $\bari A1 \ge 2$.  If the average value of $\ricci$ is positive, then $\eulgr P > 0$.
\ecoro

In a polyhedral surface we have $\bari A1 = 2$, and hence the following holds.

\coro\label{coroBAJA}
Let $K$ be a compact connected polyhedral surface.  If the average value of $\ricci$ is positive, then $\chi(K) > 0$, and hence $\pi_1 (K)$ is finite.
\ecoro

The observation in Corollary~\ref{coroBAJA} that $\pi_1 (K)$ is finite is trivial, because the classification of compact connected surfaces implies that if a compact connected polyhedral surface has positive Euler characteristic, then it is either $S^2$ or $P^2$, and in both cases it has finite fundamental group.  We stated this trivial conclusion in the corollary only for comparison with Theorem~6.1 of \cite{FORM6}.  The latter is a much stronger result, because it allows for complexes that are not surfaces, though of course the proof in \cite{FORM6} is much more substantial.

Next, we turn to the question of everywhere negative curvature on polyhedral surfaces.  For flexibility, we use polyhedral maps on surfaces, observing that any triangulation of a compact connected surface can be thought of as a polyhedral map.  If $K$ is a polyhedral map on a surface, we let $\ksim Ki$ denote the collection of all \ncells i\ of $K$, for each $i \in \{0, 1, 2\}$.

\thm\label{thmBAJ}
Let $K$ be a polyhedral map on a compact connected surface.
\enum
\item\label{thmBAJ3}
$R_0$, $\ricci$ and $R_2$ have negative values at all cells of $K$ if and only if $B_2(\eta) \ge 7$ for all $\eta \in \ksim K2$.
\item\label{thmBAJ4}
If $\chi (K) \ge 0$, then not all three of $R_0$, $\ricci$ and $R_2$ can have negative values at all cells of $K$.
\item\label{thmBAJ5}
If $\chi (K) < 0$ and $K$ is orientable, there is polyhedral map on the underlying space of $K$ such that $R_0$, $\ricci$ and $R_2$ have negative values at all cells of the polyhedral map.
\eenum
\ethm

\demo
By Lemma~\ref{lemBEBX}, we will replace $\ricci$ with $R_1$.  We ask when each of $R_0$, $R_1$ and $R_2$ has negative values.

Let $v \in \ksim K0$.  Then $R_0(v) < 0$ means $1 + \frac 32A_0(v) - [A_0(v)]^2 < 0$, and it is straightforward to verify that this condition holds if and only if $A_0(v) \ge 2$, which is always true for a polyhedral map of a surface.  Hence $R_0(v) < 0$ is always true.

Let $e \in \ksim K1$.  Because $K$ is a polyhedral map of a surface then $A_1(e) = 2$ and $B_1(e) = 2$.  Hence $R_1(e) < 0$ if and only if $16 - U_1(e) - D_1(e) < 0$, which is equivalent to $\sum_{v \covered e} A_0(v) + \sum_{\sigma \covers e} B_2(\sigma) > 16$.  The edge $e$ has two vertices, and $A_0(v) \ge 2$ for all vertices $v \in \ksim K0$, and therefore if $\sum_{\sigma \covers e} B_2(\sigma) > 12$ then $R_1(e) < 0$; this condition on $\sum_{\sigma \covers e} B_2(\sigma)$ is not necessary for obtaining $R_1(e) < 0$.  In particular, because $e$ is contained in two faces, if $B_2(\sigma) \ge 7$ for all $\sigma \in \ksim K2$ then $R_1(e) < 0$.

Let $\sigma \in \ksim K2$.  Then $R_2(\sigma) < 0$ means $1 + 6B_2(\sigma) - [B_2(\sigma)]^2 < 0$, and it is straightforward to verify that this condition holds if and only if $B_2(\sigma) \ge 7$.

Putting the above considerations together immediately implies Part~(\ref{thmBAJ3}) of this theorem.

Suppose $\chi (K) \ge 0$.  It can be verified that $K$ must have a face with no more than $6$ edges; details are left to the reader.  Part~(\ref{thmBAJ4}) of this theorem then follows immediately from Part~(\ref{thmBAJ3}).

Now suppose that $K$ is orientable and $\chi (K) < 0$.  Then the genus of $K$ is greater than or equal to $2$.  By Theorem~2(a) of \cite{M-S-W1} there is a polyhedral surface $L$ with the same genus as $K$, and with each face having $3$ edges, and each vertex contained in $7$ edges; this surface can be rectilinearly embedded in $\rrr{3}$.  (Such a polyhedral surface is called equivelar.)  Let $M$ be the dual map of $L$.  Then each face of $M$ has $7$ edges.  It follows immediately from Part~(\ref{thmBAJ3}) of this theorem that $R_0$, $R_1$ and $R_2$ have negative values at all cells of $M$, which is Part~(\ref{thmBAJ5}) of the theorem.  
\edemo

The glaring omission in Theorem~\ref{thmBAJ} is that Part~(\ref{thmBAJ5}) of the theorem treats only the orientable case.  It does not appear to be known whether for every non-orientable compact connected surface $K$ with $\chi (K) < 0$, there is a polyhedral map $M$ of the underlying space of $K$ such that every face of $M$ has at least $7$ edges.  It would be interesting to know if that holds.

\section{Comparison with D.\ Stone's Approach}
\label{secDSTO}

The approach taken in \cite{FORM6} is not the only combinatorial approach to curvature of CW complexes.  Another (in fact earlier) approach was taken in \cite{STON1} and \cite{STON2}, where an analog of Myers' Theorem was proved in dimension $2$, using ideas similar to those in \cite{FORM6}; the results of the latter are stronger than the former.  We now offer a brief comparison of our approach in the case of polyhedral surfaces with the approach of \cite{STON1} and \cite{STON2}.

In \cite{FORM6}, the \crc\ is located at the edges; in our approach the \dsc\ for polyhedral surfaces is located at all cells, though from the point of view of the analog of Myers' Theorem, the curvature of interest is located at the edges.  By contrast, the curvature defined in \cite{STON1} and \cite{STON2} is located at the vertices, and is defined as follows.  Let $K$ be a polyhedral surface, and let $v \in \ksim K0$.  Using our notation, the curvature at $v$ is given by the formula
\begin{equation}
\label{eqSTO}
R^*(v) =  2 - \sum_{\sigma > v} \left(1 - \frac 2{B_2(\sigma)}\right),
\end{equation}
where the summation is over all \ncells 2\ of $K$ that contain $v$.  

It is simple to verify that an analog of the Gauss Bonnet Theorem, which is $\sum_{v \in \ksim K0} R^*(v) = 2 \chi (K)$, holds for this definition of curvature for all polyhedral surfaces.  However, in contrast to the analog of the Gauss-Bonnet Theorem in Theorem~\ref{thmBAK}, which holds for all finite ranked posets of rank $2$, the analog of the Gauss Bonnet Theorem for $R^*$ does not hold for all \nd{2}dimensional simplicial complexes, as examples show; we omit the details.

When $K$ is a polyhedral surface, the number of \ncells 2\ that contain $v$ equals the number of edges that contain $v$, which is $A_0(v)$.  Hence we can trivially rewrite Equation~\ref{eqSTO} as 
\begin{equation}
\label{eqSTY}
R^*(v)= 2 - A_0(v) + \sum_{\sigma > v}\frac 2{B_2(\sigma)}.
\end{equation}
In contrast to the formula for $R^*$ given in Equation~\ref{eqSTO}, which is precisely as given in \cite{STON1} and \cite{STON2}, the formula given in Equation~\ref{eqSTY} does satisfy the analog of the Gauss Bonnet Theorem for all \nd{2}dimensional simplicial complexes.  However, while the formula in Equation~\ref{eqSTY} can, similarly to \crc, be extended unchanged to the context of \cfin\ ranked posets, the analog of the Gauss-Bonnet Theorem does not hold for this definition of $R^*$ for all finite ranked poset of rank $2$, as examples show; again, we omit the details.

We now return to the case where $K$ is a polyhedral surface.  The analog of the Gauss-Bonnet Theorem for $R^*$, which in contrast to Corollary~\ref{coroBBA} does not have $-1$ coefficients for any terms, implies that if $\chi (K) \ge 0$ then it cannot be the case that $R^*(w) < 0$ for all $w \in \ksim K0$, which is analogous to Theorem~\ref{thmBAJ} (\ref{thmBAJ4}).  

Similarly to one direction of Theorem~\ref{thmBAJ} (\ref{thmBAJ3}), if $B_2(\eta) \ge 7$ for all $\eta \in \ksim K2$, then
\[
R^*(v)  \le 2 - A_0(v) + \sum_{\sigma > v}\frac 27 = 2 - A_0(v) + \frac 27 A_0(v) 
\le 2 - \frac 57 \cdot 3 < 0
\]
for all $v \in \ksim K0$, because $A_0(v) \ge 3$.   (We note, however, that whereas this condition is also a necessary condition in Theorem~\ref{thmBAJ} (\ref{thmBAJ3}), it is not necessary for $R^*(v)$ for all $v \in \ksim K0$; for example, by \cite{M-S-W1}*{Theorem~2(c)} there is a polyhedral surface $M$ such that $\chi (M) = -8$, that each face has $5$ edges, and each vertex is contained in $4$ edges; it is seen that $R^*(v) < 0$ for all $v \in M$.)   

The proof of Theorem~\ref{thmBAJ} (\ref{thmBAJ5}) also shows that if $K$ is orientable and $\chi (K) < 0$, there is polyhedral surface $M$ with the same underlying space as $K$ such that $R^*(w) < 0$ for all $w \in \ksim K0$.

Finally, we note that $R^*$ does not work any better than $\ricci$, or $R_0$, $R_1$ and $R_2$ in relation to a possible analog of Myers' Theorem for posets.  Let $P$ be the poset shown in Figure~\ref{figDIRIm}, and $v \in P$ be the element shown in the figure.  The reader can verify that $R^*(v) = 3$.  It is therefore not the case that $R^*(w)$ is positive and bounded away from zero for all $w \in \rkiset P0$ guarantees that the poset is finite.


\begin{bibdiv}

\begin{biblist}[\normalsize]


\bib{A-M-M}{article}{
   author={Alsing, Paul M.},
   author={McDonald, Jonathan R.},
   author={Miller, Warner A.},
   title={The simplicial Ricci tensor},
   journal={Classical Quantum Gravity},
   volume={28},
   date={2011},
   number={15},
   pages={155007, 17}
}

\bib{BA1}{article}{
author = {Banchoff, Thomas},
title = {Critical points and curvature for embedded polyhedra},
journal = {J. Diff. Geom.},
volume = {1},
date = {1967},
pages = {245--256}
}

\bib{BJOR2}{article}{
   author={Bj{\"o}rner, A.},
   title={Posets, regular CW complexes and Bruhat order},
   journal={European J. Combin.},
   volume={5},
   date={1984},
   number={1},
   pages={7--16}
}

\bib{BL4}{article}{
author = {Bloch, Ethan D.},
title = {The angle defect for arbitrary polyhedra},
journal = {Beitr{\"a}ge Algebra Geom.},
volume = {39},
date = {1998},
pages = {379--393}
}

\bib{BR-SC}{article}{
   author={Brehm, Ulrich},
   author={Schulte, Egon},
   title={Polyhedral maps},
   conference={
      title={Handbook of discrete and computational geometry},
   },
   book={
      series={CRC Press Ser. Discrete Math. Appl.},
      publisher={CRC},
      place={Boca Raton, FL},
   },
   date={1997},
   pages={345--358}
}

\bib{BROO1}{article}{
   author={Brooks, Robert},
   title={A construction of metrics of negative Ricci curvature},
   journal={J. Differential Geom.},
   volume={29},
   date={1989},
   number={1},
   pages={85--94}
}

\bib{C-M-S}{article}{
author = {Cheeger, J.},
author = {Muller, W.},
author = {Schrader, R.},
title = {On the curvature of piecewise flat spaces},
journal = {Commun. Math. Phys.},
volume = {92},
date = {1984},
pages = {405--454}
}

\bib{FE}{book}{
author = {Federico, P. J.},
title = {Descartes on Polyhedra},
publisher = {Springer-Verlag},
address = {New York},
date = {1982}
}

\bib{FORM6}{article}{
author = {Forman, Robin},
title = {{B}ochner's method for cell complexes and combinatorial {R}icci curvature},
journal = {Discrete Comput. Geom.},
volume = {29},
number = {3},
date = {2003},
pages = {323--374}
}

\bib{FRIT}{article}{
   author={Fritz, H.},
   title={Isoparametric finite element approximation of Ricci curvature},
   journal={IMA J. Numer. Anal.},
   volume={33},
   date={2013},
   number={4},
   pages={1265--1290}
}

\bib{GAO1}{article}{
   author={Gao, L. Zhiyong},
   title={The construction of negatively Ricci curved manifolds},
   journal={Math. Ann.},
   volume={271},
   date={1985},
   number={2},
   pages={185--208}
}

\bib{GA-YA}{article}{
   author={Gao, L. Zhiyong},
   author={Yau, S.-T.},
   title={The existence of negatively Ricci curved metrics on three-manifolds},
   journal={Invent. Math.},
   volume={85},
   date={1986},
   number={3},
   pages={637--652}
}

\bib{GR2}{article}{
author = {Gr{\"u}nbaum, Branko},
title = {{G}rassmann angles of convex polytopes},
journal = {Acta. Math.},
volume = {121},
date = {1968},
pages = {293--302}
}

\bib{G-S2}{article}{
author = {Gr\"unbaum, Branko},
author = {Shephard, G. C.},
title = {{D}escartes' theorem in $n$ dimensions},
journal = {Enseign. Math. (2)},
volume = {37},
date = {1991},
pages = {11--15}
}

\bib{HIGU}{article}{
   author={Higuchi, Yusuke},
   title={Combinatorial curvature for planar graphs},
   journal={J. Graph Theory},
   volume={38},
   date={2001},
   number={4},
   pages={220--229}
}

\bib{LOHK1}{article}{
   author={Lohkamp, Joachim},
   title={Negative bending of open manifolds},
   journal={J. Differential Geom.},
   volume={40},
   date={1994},
   number={3},
   pages={461--474}
}

\bib{LOHK2}{article}{
   author={Lohkamp, Joachim},
   title={Metrics of negative Ricci curvature},
   journal={Ann. of Math. (2)},
   volume={140},
   date={1994},
   number={3},
   pages={655--683}
}

\bib{L-S}{article}{
   author={Luo, Feng},
   author={Stong, Richard},
   title={Combinatorics of triangulations of $3$-manifolds},
   journal={Trans. Amer. Math. Soc.},
   volume={337},
   date={1993},
   number={2},
   pages={891--906}
}

\bib{M-S-W1}{article}{
   author={McMullen, P.},
   author={Schulz, Ch.},
   author={Wills, J. M.},
   title={Equivelar polyhedral manifolds in $E\sp{3}$},
   journal={Israel J. Math.},
   volume={41},
   date={1982},
   number={4},
   pages={331--346}
}

\bib{MYER}{article}{
   author={Myers, S. B.},
   title={Riemannian manifolds with positive mean curvature},
   journal={Duke Math. J.},
   volume={8},
   date={1941},
   pages={401--404}
}

\bib{ROS00}{collection}{
   title={Handbook of discrete and combinatorial mathematics},
   editor={Rosen, Kenneth H.},
   editor={Michaels, John G.},
   editor={Gross, Jonathan L.},
   editor={Grossman, Jerrold W.},
   editor={Shier, Douglas R.},
   publisher={CRC Press},
   place={Boca Raton, FL},
   date={2000}
}

\bib{SH2}{article}{
author = {Shephard, G. C.},
title = {Angle deficiencies of convex polytopes},
journal = {J. London Math. Soc.},
volume = {43},
date = {1968},
pages = {325--336}
}

\bib{STAN}{book}{
   author={Stanley, Richard P.},
   title={Enumerative combinatorics. Vol. 1},
   series={Cambridge Studies in Advanced Mathematics},
   volume={49},
   note={With a foreword by Gian-Carlo Rota;
   Corrected reprint of the 1986 original},
   publisher={Cambridge University Press},
   place={Cambridge},
   date={1997}
}
	
\bib{STON1}{article}{
   author={Stone, David A.},
   title={A combinatorial analogue of a theorem of Myers},
   journal={Illinois J. Math.},
   volume={20},
   date={1976},
   number={1},
   pages={12--21}
}

\bib{STON2}{article}{
   author={Stone, David A.},
   title={Correction to my paper: ``A combinatorial analogue of a theorem of
   Myers'' (Illinois J. Math. {\bf 20} (1976), no. 1, 12--21)},
   journal={Illinois J. Math.},
   volume={20},
   date={1976},
   number={3},
   pages={551--554}
}	


\end{biblist}

\end{bibdiv}
\end{document}